\documentclass[11pt]{amsart}
\usepackage{amssymb}
\usepackage{latexsym,bm}

\numberwithin{equation}{section} \theoremstyle{section}
\newtheorem{De}[equation]{Definition}

\theoremstyle{plain}
\newtheorem{Example}[equation]{Example}
\newtheorem{Prop}[equation]{Proposition}
\newtheorem{Theorem}[equation]{Theorem}

\newtheorem{Lemma}[equation]{Lemma}
\newtheorem{Cor}[equation]{Corollary}

\title{The Tracial Rokhlin Property for  Automorphisms on Non-Simple
$C$*-algebras}

\author{JIAJIE HUA}

\date{}

\begin{document}

\maketitle \markboth{JIAJIE HUA}{THE TRACIAL ROKHLIN PROPERTY FOR
AUTOMORPHISMS ON $C$*-ALGEBRAS}
\renewcommand{\thefootnote}{{}}
\footnote{\hspace{-14pt} {\it 2000 Mathematics Subject Classification}. Primary 46L55: Secondary 46L35, 46L40.This research is part of the author¡¯s Ph.D.
thesis at East China Normal University, completed under the
direction of Huaxin Lin. The author was supported the National
Natural Science Foundation of China (Nos. 10771069, 10671068) and
Shanghai Priority Academic Discipline Programme (No.
  B407).}
\begin{abstract} Let $A$ be a unital AF-algebra (simple or non-simple) and let $\alpha$ be an automorphism of $A$. Suppose that $\alpha$ has
certain Rokhlin property and $A$ is $\alpha$-simple. Suppose also
that there is an integer $J\geq1$ such that
$\alpha^{J}_{*0}=$id$_{K_{0}(A)}$, we show that
$A\rtimes_{\alpha}\mathbb{Z}$ has tracial rank zero.
\end{abstract}

\section{Introduction} We introduce  certain Rokhlin property for
automorphisms on unital $C$*-algebras. The Rokhlin property in
ergodic theory was adopted to the context of von Neumann algebras by
Connes \cite{A.Connes}. It was adopted by Herman and Oeneanu
\cite{R.Herman} for UHF-algebras. R$\phi$rdam \cite{M.R} and
Kishimoto \cite{A.Kishimoto2} introduced the Rokhlin property to a
much more general context of $C$*-algebras, then Osaka and Phillips
studied integer group actions which satisfy certain type of Rokhlin
property on some simple $C$*-algebras \cite{Osaka}. More recently,
Lin studied the Rokhlin property for automorphisms on simple
$C$*-algebras \cite{H.Lin3}.

Phillips proposed that how to introduce appropriate  Rokhlin
property to non-simple $C$*-algebras. In this paper we attempt to
introduce certain Rokhlin property to  non-simple $C$*-algebras,
when $C$*-algebra is simple, this Rokhlin property is weaker than
the Rokhlin property in \cite{H.Lin3,Osaka}. If an integer group
action of  a $C$*-algebra has this Rokhlin property, we can conclude
that its crossed product is in the $C$*-algebra class of tracial
rank zero. In particular, these algebras all belong to the class
known currently to be classifiable by K-theoretic invariants in the sense of
the Elliott classification program.  We
hope that this case will lead us to more interesting in the Rokhlin
property to non-simple $C$*-algebras.

The organization of the paper is as follows. In Section 1, we
briefly recall the notion of $C$*-algebras, then we introduce
certain Rokhlin property and discuss some property of crossed
product $A\rtimes_{\alpha}\mathbb{Z}$ when an automorphism $\alpha$
of a $C$*-algebra $A$ has the Rokhlin property. In Section 2, we
show that if $A$ is a unital AF-algebra, suppose that
$\alpha\in$Aut$(A)$ has the tracial cyclic Rokhlin property and $A$
is $\alpha$-simple, suppose also that there is an integer $J\geq1$
such that $\alpha^{J}_{*0}=$id$_{K_{0}(A)}$. Then
$A\rtimes_{\alpha}\mathbb{Z}$ has tracial rank zero.

\section{The Tracial Rokhlin Property}

We will use the following convention:\\

(1) Let $A$ be a $C$*-algebra, let $a\in A$ be a positive element
and let $p\in A$ be a projection. We write $[p]\leq [a]$ if there is
a projection $q\in \overline{aAa}$ and a partial isometry $v\in A$
such that $v^{*}v=p$ and $vv^{*}=q.$

(2) Let $A$ be a $C$*-algebra. We denote by Aut$(A)$ the
automorphism group of $A$. If $A$ is unital and $u\in A$ is a
unitary, we denote by ad$u$ the inner automorphism defined by
ad$u(a)=u^{*}au$ for all $a\in A.$

(3) Let $x\in A,\varepsilon>0$ and $\mathcal{F}\subset A.$ We write
$x\in_{\varepsilon}\mathcal{F},$ if
dist$(x,\mathcal{F})<\varepsilon$, or there is $y\in \mathcal{F}$
such that $\|x-y\|<\varepsilon.$

(4) Let $A$ be a $C$*-algebra and $\alpha\in$ Aut$(A)$. We say $A$
is $\alpha$-simple if $A$ does not have any non-trivial
$\alpha$-invariant closed two-sided ideals.

(5) A unital $C$*-algebra is said to have real rank zero, written
RR$(A)=0$, if the set of invertible self-adjoint elements is dense
in self-adjoint elements of $A.$ Note that every unital AF-algebra
has real rank zero.

(6) A unital $C$*-algebra $A$ has the (SP)-property if every
non-zero hereditary $C$*-subalgebra of $A$ has a non-zero
projection. Note that every $C$*-algebra $A$ with real rank zero has
the (SP)-property.

(7) Let $T(A)$ be the tracial state space of a unital $C$*-algebra
$A$. It is a compact convex set.

(8) we say the order on projection over a unital $C$*-algebra $A$ is
determined by traces, if for any two projections $p,q\in A,
\tau(p)<\tau(q)$ for all $\tau\in T(A)$ implies that $p$ is
equivalent to a projection $p'\leq q.$

\begin{De}  We denote by $\mathcal{I}^{(0)}$ the class of all finite
dimensional $C$*-algebras, and denote by $\mathcal{I}^{(k)}$ the
class of all unital $C$*-algebras which are unital hereditary
$C$*-subalgebras of $C$*-algebras of the form $C(X)\otimes F$, where
$X$ is a $k$-dimensional finite CW complex and $F\in
\mathcal{I}^{(0)}.$
\end{De}

 We recall the definition of tracial topological rank of
$C$*-algebras.

\begin{De} {\rm([8])} Let $A$ be a unital simple $C$*-algebra.
Then $A$ is said to have tracial (topological) rank no more than $k$
if for any $\varepsilon>0,$ any finite set $\mathcal{F}\subset
A,$ and any non-zero positive element $a\in A,$ there exist a
nonzero projection $p\in A$ and a $C$*-subalgebra $B\in
\mathcal{I}^{(k)}$ with $1_B=p$ such that:

(1) $\|px-xp\|<\varepsilon$ for all $x\in \mathcal{F}$.

(2) $pxp\in_{\varepsilon} B$ for all $x\in \mathcal{F}$.

(3) $[1-p]\leq [a].$
\end{De}

If $A$ has tracial rank no more than $k$, we will write
$\mathrm{TR}(A)\leq k.$ If furthermore, $\mathrm{TR}(A)\nleq k-1,$
then we say $\mathrm{TR}(A)=k.$

\begin{De}  Let $A$ be a unital $C$*-algebra and let
$\alpha\in$Aut$(A)$. Let $a\in A$ be a positive element and let
$p\in A$ be a projection.  We say $[p]\leq_{\alpha} [a]$ if there
exist the mutually orthogonal projections $p_{i},$ the mutually
orthogonal positive elements $a_{i}$ and $s_{i}\in \mathbb{Z}$ for
$i=1,2,\cdots,n$ such that
$p=\sum_{i=1}^{n}p_{i},\{a_{i}\}_{i=1}^{n}$ belong to the hereditary
$C$*-subalgebra generated by $a$,  and $[\alpha^{s_{i}}(p_{i})]\leq
[a_{i}],$ $i=1\cdots n.$
\end{De}

 By this definition, we can compare non-zero positive elements with full positive
 elements by the action of  $\alpha$.

\begin{Example} Let  $A=A_{0}\oplus A_{0}$, where $A_{0}$ is  an infinite dimensional unital simple $C$*-algebra with real rank zero,
 let $\alpha\in$Aut$(A)$ such that
$\alpha(a_{0},b_{0})=(b_{0},a_{0})$, where $a_{0},b_{0}\in A_{0}$,
then for any non-zero projection $q\in A$, there exists a projection
$p=(p_{1},p_{2})\in A ,p_{1}\neq 0,p_{2}\neq 0$ such that
$[p]\leq_{\alpha}[q].$ \end{Example}

\begin{De}  Let $A$ be a unital $C$*-algebra and let
$\alpha\in$Aut$(A)$. We say $\alpha$ has the tracial Rokhlin
property if for every $\varepsilon>0$, every $n\in \mathbb{N}$,
every nonzero positive element $a\in A$, every finite set
$\mathcal{F}\subset A, \mathcal{F}=\{p_{1},\cdots,p_{m},a_{1},$
$\cdots,a_{s}\},$ where $\{p_{i}\},i=1,\cdots,m$ are the mutually
orthogonal projections, there are the mutually orthogonal
projections $e_{1},e_{2},\cdots,e_{n}\in A$ such that:

 (1) $\|\alpha(e_{j})-e_{j+1}\|<\varepsilon$ for $1\leq j \leq
n-1.$

(2) $\|e_{j}b-be_{j}\|<\varepsilon$ for $1\leq j \leq n$ and all
$b\in \mathcal{F}.$

(3) $\|e_{1}p_{j}e_{1}\|\geq 1-\varepsilon$ for $1\leq j \leq m.$

(4) With $e=\sum_{j=1}^{n}e_{j},[1-e]\leq_{\alpha}[a].$
\end{De}

When $A$ is a unital simple $C$*-algebra, above definition is weaker
than the Rokhlin property in \cite{H.Lin3,Osaka}. We weak the
condition (4) to require positive elements can be compared by action
of $\alpha$.

 We define a slightly stronger version of the tracial Rokhlin
 property.
\begin{De}  Let $A$ be a unital $C$*-algebra and let
$\alpha\in$Aut$(A)$. We say $\alpha$ has the  tracial cyclic Rokhlin
property if for every $\varepsilon>0$, every $n\in \mathbb{N}$,
every nonzero positive element $a\in A$, every finite set
$\mathcal{F}\subset A, \mathcal{F}=\{p_{1},\cdots,p_{m}, a_{1},$
$\cdots,a_{s}\},$ where $\{p_{i}\},i=1,\cdots,m$ are the mutually
orthogonal projections, there are the mutually orthogonal
projections $e_{1},e_{2},\cdots,e_{n}\in A$ such that:

(1) $\|\alpha(e_{j})-e_{j+1}\|<\varepsilon$ for $1\leq j \leq n,$
where $e_{n+1}=e_{1}.$

(2) $\|e_{j}b-be_{j}\|<\varepsilon$ for $1\leq j \leq n$ and all
$b\in \mathcal{F}.$

(3) $\|e_{1}p_{j}e_{1}\|\geq 1-\varepsilon$ for $1\leq j \leq m.$

(4) With $e=\sum_{j=1}^{n}e_{j},[1-e]\leq_{\alpha}[a].$
\end{De}

The only difference between the tracial Rokhlin property and the
tracial cyclic Rokhlin property is that in condition (1), we require
that $\|\alpha(e_{n})-e_{1}\|<\varepsilon.$

\begin{Theorem}  Let $A$ be a unital $C$*-algebra with real rank zero,  let
$\alpha\in$Aut$(A)$ have the tracial Rokhlin property. Then $A$ is
$\alpha$-simple if and only if the crossed product
$A\rtimes_{\alpha} \mathbb{Z}$ is simple.
\end{Theorem}

\begin{proof} Let $I$ be an $\alpha$-invariant norm closed two-sided
ideal of $A$. Then $I\rtimes_{\alpha}\mathbb{Z}$ is a norm closed
two-sided ideal of $A\rtimes_{\alpha} \mathbb{Z}$ by Lemma 1 of
\cite{Jang}.

Conversely, let $a$ be a positive element of the $C$*-algebra $A$,
$\mathcal{F}=\{a_{i};i=1,2,\cdots.n\}$ elements of $A$, $s_{i}\in
\mathbb{N}, i=1,2,\cdots,n$ and $\varepsilon>0$, we prove that there
exists a positive element $x\in A$ with $\|x\|=1$ such that
$$\|xax\|\geq\|a\|-\varepsilon, \quad
\|xa_{i}\alpha^{s_{i}}(x)\|\leq \varepsilon, \quad
\|xa_i-a_ix\|<\varepsilon, i=1,2,\cdots,n.\quad\quad (*)$$

Because $A$ has real rank zero, let $\varepsilon>0,$ by Theorem
3.2.5 of \cite{H.Lin2}, there are mutually orthogonal projections
$p_{1},p_{2},\cdots,p_{m}$ and positive real numbers
$\lambda_{1},\lambda_{2},\cdots,\lambda_{m}$ such that
$\|a-\sum_{i=1}^{m}\lambda_{i}p_{i}\|<\varepsilon/3,$ let
$a_{0}=\sum_{i=1}^{m}\lambda_{i}p_{i}$, $C=\max\{\|a_i\|
|i=1,\dots,n\}$ and $N=\max\{s_{1},s_{2},\cdots,s_{n}\}.$

 Let
$\varepsilon_{0}=\min\{\frac{\varepsilon}{3\|a_{0}\|},
\frac{\varepsilon}{(N+2)C}\}.$

Apply the tracial Rokhlin property with  $N$ in place of n, with
$\varepsilon_{0}$ in place of $\varepsilon.$ We can obtain
$e_{1},e_{2}, \cdots,e_{N},$ such that

(1) $\|\alpha(e_{j})-e_{j+1}\|<\varepsilon_{0}$ for $1\leq j \leq
N-1,$

(2) $\|e_{j}a_{i}-a_{i}e_{j}\|<\varepsilon_{0}$ for $1\leq j \leq N$
and $1\leq i\leq n$,

(3) $\|e_{1}p_{j}e_{1}\|\geq 1-\varepsilon_{0}$ for $1\leq j \leq
m,$

 then\\
$\|e_{1}a_{0}e_{1}\|=\|\sum_{i=1}^{m}\lambda_{i}e_{1}p_{i}e_{1}\|\geq
\|\lambda_{i}e_{1}p_{i}e_{1}\|\geq\lambda_{i}(1-\varepsilon_{0}),i=1,2,\cdots,m$.

we get
$$\|e_{1}a_{0}e_{1}\|\geq\|a_{0}\|(1-\varepsilon_{0})\geq\|a_{0}\|-\frac{\varepsilon}{3}.$$

then
\begin{eqnarray*}\|e_{1}ae_{1}\|&=&\|e_{1}a_{0}e_{1}+e_{1}ae_{1}-e_{1}a_0e_{1}\|\geq
\|e_{1}a_{0}e_{1}\|-\|e_{1}ae_{1}-e_{1}a_{0}e_{1}\|\\
&\geq& \|e_{1}a_{0}e_{1}\|-\frac{\varepsilon}{3}\geq
\|a_{0}\|-\frac{\varepsilon}{3}-\frac{\varepsilon}{3}\geq
\|a\|-\frac{\varepsilon}{3}-\frac{\varepsilon}{3}-\frac{\varepsilon}{3}=\|a\|-\varepsilon.
\end{eqnarray*}
\begin{eqnarray*}& &\|e_{1}a_{i}\alpha^{s_{i}}(e_{1})\|\\&=&\|e_{1}a_{i}\alpha^{s_{i}}(e_{1})-e_{1}a_{i}\alpha^{s_{i}-1}(e_{1})+e_{1}a_{i}\alpha^{s_{i}-1}(e_{1})
+\cdots+e_{1}a_{i}\alpha^{1}(e_{1})\|\\
&<&\|e_{1}a_{i}\alpha^{1}(e_{1})\|+(s_{i}-1)\varepsilon_{0}\|a_i\|<
\|a_{i}e_{1}\alpha^{1}(e_{1})\|+s_{i}\varepsilon_{0}\|a_i\|\\
&<& (s_{i}+1)\varepsilon_{0}\|a_i\|<\varepsilon.\end{eqnarray*}

So we get ($*$). Apply this condition and $A$ is $\alpha$-simple, we
can complete the proof as same as Theorem 3.1 of
\cite{A.Kishimoto1}, we omit it.
\end{proof}

 Apply the $(*)$ and the same proof of Theorem 4.2 of \cite{J.Jeong1}, we can get the following
 result.
\begin{Theorem} Let $A$ be a unital $C$*-algebra with real rank
zero and let $\alpha\in$Aut$(A)$ have the tracial Rokhlin property
and $A$ is $\alpha$-simple. Then any non-zero hereditary
$C$*-subalgebra of the  crossed product
$A\rtimes_{\alpha}\mathbb{Z}$ has a non-zero projection which is
equivalent to a projection in $A.$
\end{Theorem}

\begin{Lemma} Let $B=M_{r(1)}\oplus M_{r(2)}\oplus \cdots \oplus M_{r(l)}$ be a finite dimensional $C$*-subalgebra of a
unital $C$*-algebra $A$, Let $e^{(s)}_{i,j}\in B$ be a system of
matrix units for $M_{r(s)},s=1,2,\cdots, l.$ Then for any
$\delta>0,$ there exists $\sigma>0$ satisfying the following: If
$\|pe^{(s)}_{i,i}-e^{(s)}_{i,i}p\|<\sigma$ and
$\|pe^{(s)}_{i,i}p\|>1/2$ for $s=1,2,\dots,l, i=1,2,\dots,r(s),$
then there is an monomorphism $\varphi: B\rightarrow pAp$ such that
$\|pbp-\varphi(b)\|<\delta \|b\|$ for all $b\in B.$
\end{Lemma}

\begin{proof} It follows from the argument in section 2.5 of
\cite{H.Lin2} and Proposition 2.3 of \cite{H.Lin4}.
\end{proof}

\begin{Prop} Let $A$ be a unital $C$*-algebra, Suppose that
$\alpha\in$Aut$(A)$ is approximately inner and has the tracial
Rokhlin property, if for any closed two-sided  ideal $I$ of
$C$*-algebra $A$, there is $n\in \mathbb{N}$, $n$ only depends on
$I$,  such that $K_{0}(A/I)$ is not $n$-divisible, then $A$ is
$\alpha$-simple.
\end{Prop}

\begin{proof} Suppose that $A$ is not $\alpha$-simple, so there
exists a closed two-sided  ideal $I$ of $C$*-algebra $A$ such that
$\alpha(I)=I.$ By the hypothesis, there is $n\in \mathbb{N}$ such
that $K_{0}(A/I)$ is not $n$-divisible.

Let $a\in I$ be a non-zero positive element, and $0<\varepsilon<1$,
there are the mutually orthogonal projections
$e_{1},e_{2},\cdots,e_{n}\in A$ such that

(1) $\|\alpha(e_{j})-e_{j+1}\|<\varepsilon$ for $1\leq j \leq n-1,$

(2) With $e=\sum_{j=1}^{n}e_{j},[1-e]\leq_{\alpha}[a].$

Because $\alpha$ is approximately inner and (1), we have
$[e_{1}]=[e_{2}]=\cdots=[e_{n}]$ in $K_{0}(A).$

If $p\in A$ is a projection such that $[p]\leq[b],$ where $b\in I$
is a positive element, then there is a $v\in A$ such that $v^{*}v=p$
and $vv^{*}\in \overline{bAb}\subset I$,  if $\pi: A\rightarrow A/I$
denotes quotient map,  $\pi(v)\pi(v^{*})=0$ in $A/I,$ $\pi(v)=0$ in
$A/I,$ then $p\in I.$

In (2), $[1-e]\leq_{\alpha}[a]$, by the definition of
$\leq_{\alpha}$, $a\in I$ and the discuss above, we have $1-e\in I,$
so $\pi(1-e)=0,$ $[1-e]=0$ in $K_{0}(A/I),$ then $n[e_{1}]=[1]$ in
$K_{0}(A/I),$ this is contradictory to $K_{0}(A/I)$ is not
$n$-divisible.

\end{proof}

\section{Main result}

In the proof of Theorem 3.3, we first prove
$TR(A\rtimes_{\alpha}\mathbb{Z})\leq1,$ then use the following Lemma
3.1 to prove $RR(A\rtimes_{\alpha}\mathbb{Z})=0$. The following
Lemma is similar to Lemma 2.5 of \cite{Osaka}.

\begin{Lemma} Let $A$ be a unital $C$*-algebra with real rank zero
 and let $\alpha\in$Aut$(A)$ have the tracial Rokhlin property.
Suppose that $A$ is $\alpha$-simple and the order on projection over
$A\rtimes_{\alpha}\mathbb{Z}$ is determined by traces. Let $\iota:
A\rightarrow A\rtimes_{\alpha}\mathbb{Z}$ be the inclusion map. Then
for every finite set $F\subset A\rtimes_{\alpha}\mathbb{Z},$ every
$\varepsilon>0,$ every nonzero positive element $z\in
A\rtimes_{\alpha}\mathbb{Z},$ and every sufficiently large $n\in
N$(depending on $F,\varepsilon$ and $z$), there exist a projection
$e\in A\subset A\rtimes_{\alpha}\mathbb{Z},$ a unital subalgebra
$D\subset e(A\rtimes_{\alpha}\mathbb{Z})e,$ a projection $p\in D$, a
projection $f\in A,$ and an isomorphism $\varphi:M_n\otimes
fAf\rightarrow D,$ such that:

(1) With $(e_{j,k})$ being the standard system of matrix units for
$M_n$, we have $\varphi(e_{1,1}\otimes a)=\iota(a)$ for all $a\in
fAf$ and $\varphi(e_{k,k}\otimes 1)\in \iota (A)$ for $1\leq k \leq
n.$

(2) With $(e_{j,k})$ as in $(1)$, we have $\|\varphi(e_{j,j}\otimes
a)-\alpha^{j-1}(\iota(a))\|\leq \varepsilon\|a\|$ for all $a\in
fAf.$

(3) For every $a\in F$, there exists $b_1,b_2\in D$ such that
$\|pa-b_1\|<\varepsilon, \|ap-b_2\|<\varepsilon,$ and
$\|b_1\|,\|b_2\|\leq \|a\|.$

(4) There is $m\in \mathbb{N}$ such that $2m/n<\varepsilon$ and
$p=\sum_{j=m+1}^{n-m}\varphi(e_{j,j}\otimes 1).$

(5) The projection $1-p$ is Murray-von Neumann equivalent in
$A\rtimes_{\alpha}\mathbb{Z}$ to a projection in the hereditary
subalgebra of $A\rtimes_{\alpha}\mathbb{Z}$ generated by $z$ and
$\tau(1-p)<\varepsilon$ for all $\tau\in
T(A\rtimes_{\alpha}\mathbb{Z}).$
\end{Lemma}

\begin{proof} Let $\varepsilon>0,$ and let $F\subset
A\rtimes_{\alpha}\mathbb{Z}$ be a finite set. Let $z\in
A\rtimes_{\alpha}\mathbb{Z}$ be a nonzero positive element.

Let $u$ be the standard unitary in the crossed product
$A\rtimes_{\alpha}\mathbb{Z}$. We regard $A$ as a subalgebra of
$A\rtimes_{\alpha}\mathbb{Z}$ in the usual way. Choose $m\in
\mathbb{N}$ such that for every $x\in F$ there are $a_l\in A$ for
$-m\leq l\leq m$ such that $\|x-\sum_{l=-m}^m a_l
u^l\|<\frac{\varepsilon}{2}.$ For each $x\in F$ choose one such
expression, and let $S\subset A$ be a finite set which contains all
the coefficients used for all elements of $F$. Let $M=1+\sup_{a\in
S}\|a\|.$

Since $A\rtimes_{\alpha}\mathbb{Z}$ has (SP)-property and is simple
by Theorem 2.8 and Theorem 2.7, we can apply Lemma 3.5.7 of
\cite{H.Lin2} to find nonzero orthogonal Murray-von Neumann
equivalent projections $g_0,g_1,\dots, g_{2m} \in
z(A\rtimes_{\alpha}\mathbb{Z})z.$

Since $A\rtimes_{\alpha}\mathbb{Z}$ is simple, $g_0$ is a nonzero
projection, and the tracial state space
$T(A\rtimes_{\alpha}\mathbb{Z})$ of $A\rtimes_{\alpha}\mathbb{Z}$ is
weak-$*$ compact, we have $$\delta=\inf_{\tau\in
T(A\rtimes_{\alpha}\mathbb{Z})}\tau(g_0)>0.$$ Now let $n\in
\mathbb{N}$  be any integer such that
$n>\max(\frac{1}{\delta},(N+2)(2m+1),\frac{4m}{\varepsilon}).$

Set $\varepsilon_0=\frac{\varepsilon}{10(2m+1)n^2M}.$

Choose $\varepsilon_1>0$ so small that whenever $e_1,e_2,\dots,e_n$
are mutually orthogonal projections in a unital $C$*-algebra $B$ and
$u\in B$ is a unitary such that $\|ue_ju^*-e_{j+1}\|<\varepsilon_1$
for $1\leq j\leq n,$ then there is a unitary $v\in B$ such that
$\|v-u\|<\varepsilon_0$ and $ve_jv^*=e_{j+1}$ for $1\leq j\leq n$.
We can apply Lemma 3.5.7 of \cite{H.Lin2} to find nonzero orthogonal
Murray-von Neumann equivalent projections $h_1,h_2,\dots,h_{n+2}\in
\overline{g_0(A\rtimes_{\alpha}\mathbb{Z})g_0}$
 which are Murray-von Neumann
equivalent in $A\rtimes_{\alpha}\mathbb{Z}.$ Further apply Theorem
2.8 to find a nonzero projection $q\in A$ which is Murray-von
Neumann equivalent in $A\rtimes_{\alpha}\mathbb{Z}$ to a projection
in $\overline{h_1(A\rtimes_{\alpha}\mathbb{Z})h_1.}$

Apply the tracial Rokhlin property with $n-1$  in place of $n$, with
$S$ in place of $F$, with $\min(1,\varepsilon_0,\varepsilon_1)$ in
place of $\varepsilon$, and with $q$ in place of $x.$ Call the
resulting projections $e_1, e_2,\dots,e_n,$
 and let $e=\sum_{j=1}^n e_j,$ $[1-e]\leq_\alpha [q].$ Apply the choice
of $\varepsilon_1$ to these projections and the standard unitary
$u,$ obtaining a unitary $v\in A\rtimes_{\alpha}\mathbb{Z}$ as in
the previous paragraph.

We can get the conditions (1),(2),(3),(4) by the same proof of Lemma
2.5 of \cite{Osaka}. We omit them.

It remains to verify Condition (5) of the conclusion. We have
$$1-p=1-e+\sum_{j=1}^me_j+\sum_{j=n-m+1}^ne_j.$$
By constrction we have $[1-e]\leq_\alpha [h_1]\leq [g_0].$ Now let
$\tau$ be any tracial state on $A\rtimes_{\alpha}\mathbb{Z}.$ Then
$\tau(e_j)=\tau(e_1)$ for all $j$, whence $\tau(e_j)\leq
\frac{1}{n}.$ The inequality $n>\frac{1}{\delta}\geq
\frac{1}{\tau(g_0)}$ therefore implies $\tau(e_j)<\tau(g_0).$ Since
all $g_j$ are Murray-von Neumann equivalent, it follows that for any
tracial state $\tau$ on $A\rtimes_{\alpha}\mathbb{Z},$ we have
$\tau(e_j)<\tau(g_j)$ and $\tau(e_{n-j})<\tau(g_{m+j})$ for $1\leq
j\leq m.$ So the order on projection over
$A\rtimes_{\alpha}\mathbb{Z}$ is determined by traces implies that
$e_j\leq g_j$ and $e_{n-j}\leq g_{m+j}$ in
$A\rtimes_{\alpha}\mathbb{Z}$ for $1\leq j\leq m.$ Thus
$[1-p]\leq_\alpha[\sum_{j=0}^{2m}g_j]$ which is a projection in the
hereditary subalgebra $\overline{z(A\rtimes_{\alpha}\mathbb{Z})z}.$

$$\tau(1-p)=\tau(1-e)+\tau(\sum_{j=1}^me_j+\sum_{j=n-m+1}^ne_j)\leq \frac{1}{2m(n+2)}+\frac{2m}{n}<\varepsilon.$$

This is Condition  (5) of the conclusion.
\end{proof}

\begin{Theorem} Let $A$ be a unital $C$*-algebra with real rank zero
 and let $\alpha\in$Aut$(A)$ have the tracial Rokhlin property.
Suppose that $A$ is $\alpha$-simple and the order on projection over
$A\rtimes_{\alpha}\mathbb{Z}$ is determined by traces. Then
$A\rtimes_{\alpha}\mathbb{Z}$ has real rank zero.
\end{Theorem}

\begin{proof} By applying Lemma 3.1 and the same proof of Theorem
4.5 of \cite{Osaka}.
\end{proof}

\begin{Theorem} Let $A$ be a unital AF-algebra, Suppose that
$\alpha\in$Aut$(A)$ has the tracial cyclic Rokhlin property and $A$
is $\alpha$-simple. Suppose also that there is an integer $J\geq1$
such that $\alpha^{J}_{*0}=$id$_{K_{0}(A)}$. Then
TR$(A\rtimes_{\alpha}\mathbb{Z})=0.$
\end{Theorem}

\begin{proof} By Theorem 2.7, $A\rtimes_{\alpha}\mathbb{Z}$ is a
unital simple $C$*-algebra.

Let $0<\varepsilon<1$ and $\mathcal{F}\subset
A\rtimes_{\alpha}\mathbb{Z}$ be a finite set. To simplify notation,
without loss of generality, we may assume that
$\mathcal{F}=\mathcal{F}_{0}\cup \{u\},$ where
$\mathcal{F}_{0}\subset A$ is a finite subset of the unit ball which
contains $1_{A}$ and $u$ is a unitary which implements $\alpha$,
i.e., $\alpha(a)=u^{*}au$ for all $a\in A.$ Choose an integer $k$
which is a multiple of $J$ such that $2\pi/(k-2)<\varepsilon/16.$
Put $\mathcal{F}_{1}=\mathcal{F}_{0}\cup \{u^{i}a(u^{*})^{i}: a\in
\mathcal{F}_{0}, -k\leq i\leq k\}.$

Fix $b_{0}\in(A\rtimes \mathbb{Z})_{+}\backslash\{0\}.$ It follows
from Theorem 2.8 that there is a nonzero projection $r_{0}\in A$
which is equivalent to a nonzero projection in the hereditary
$C$*-subalgebra generated by $b_{0}$.

Let $\delta=\varepsilon/16k^{2}.$ Since $A$ is a unital AF-algebra,
denoted by $A=\overline{\cup_{m=1}^{\infty}A_{m}}$, where $A_m$ is a
finite-dimensional $C$*-algebra for $m=1,2,\dots,$ then there is a
lager enough $m\in\mathbb{N}$ such that $b\in_{\delta}A_{m}$ for all
$b\in \mathcal{F}_{1}$ and $1_{A}\in A_{m}.$ Let
$A_{m}=M_{r(1)}\oplus M_{r(2)}\oplus\cdots\oplus M_{r(l)}.$ Note
$[(u^{k})^{*}eu^{k}]=[e]$ in $K_{0}(A)$ for all projection $e\in
A_{m}.$ By Theorem 3.4.6 of \cite{H.Lin2}, there exists a unitary
$w\in U(A)$ such that $w^{*}(u^{k})^{*}bu^{k}w=b$ for all $b\in
A_{m}.$ Because $A$ is a AF-algebra, $w\in U_{0}(A)$. By Lemma 2.6
of \cite{H.Lin3},  we have the unitaries $w_{i},i=1,2,\cdots,k-1$
associated with finite dimensional $C$*-subalgebra $A_{m}$ such that
$w=w_{1}w_{2}\cdots w_{k-1}, \|w_{i}-1\|\leq \pi/(k-2).$ Let
$\mathcal{G}_{0}$ be a finite subset of $A_{m}$ which, for each
$b\in \mathcal{F}_{1}$ contains an element $a(b)$ such that
$\|a(b)-b\|<\delta,$ contains a systems of matrix units for each
simple summand of $A_{m}.$

Define $\mathcal{F}_{2}=\{u^{i}bu^{-i}: b\in \mathcal{G}_{0}, -k\leq
i \leq k\}$ and let $w_{k}=1$

$\mathcal{F}_{3}=\{(w_{i_{1}}w_{i_{1}+1}\cdots
w_{i})a(w_{i_{2}}w_{i_{2}+1}\cdots w_{i})^{*}: a\in
\mathcal{F}_{1}\cup \mathcal{F}_{2}, 1 \leq i,i_{1},i_{2}\leq k ,
i_{1}\leq i,i_{2}\leq i \}.$ Note that $w,w_{i}\in
\mathcal{F}_{3},i=1,2,\cdots,k-1.$

Since $\alpha$ has the tracial cyclic Rokhlin property, let
$e_{i,j}^{(s)}\in A_{m}$ be a system of matrix units for
$M_{r(s)},s=1,2,\cdots,l.$ Let $\sigma>0$ be associated with $A_{m}$
and $\delta$ in Lemma 2.9. Let $\eta=\min\{\delta,\sigma\},$ there
exist projections $e_{1},e_{2},\cdots,e_{k}\in A$ such that:

(1) $\|\alpha(e_{i})-e_{i+1}\|<\eta/k$ for $1\leq i \leq
k,e_{k+1}=e_{1}$.

(2) $\|e_{i}a-ae_{i}\|<\eta/k$ for $a\in \mathcal{F}_{3}.$

(3) $\|e_{1}e_{jj}^{(s)}e_{1}\|\geq 1-\eta/k$ for $s=1,2,\cdots,l,
j=1,2,\cdots,r(s).$

(4) $[1-\sum_{i=1}^{k}e_{i}]\leq_{\alpha}[r_{0}].$

Set $p=\sum_{i=1}^{k}e_{i}.$ From (1) above, one estimates that
\begin{eqnarray*}\|up-pu\|&=&\|\sum_{i=1}^{k}ue_{i+1}-\sum_{i=1}^{k}e_{i}u\|\leq\sum_{i=1}^{k}\|ue_{i+1}-e_{i}u\|
\\ &=&\sum_{i=1}^{k}\|ue_{i+1}-u\alpha(e_{i})\|<\eta.
\end{eqnarray*}

By (1) above, one sees that there is a unitary $v\in A$ such that
$\|v-1\|<2\eta/k$ and $v^{*}u^{*}e_{i}uv=e_{i+1}, i=1,2,\cdots,k.$
Set $u_{1}=uv.$ Then $u_{1}^{*}e_{i}u_{1}=e_{i+1},i=1,2,\cdots,k$
and $e_{k+1}=e_{1}$. In particular, $u_{1}^{k}e_{1}=e_{1}u_{1}^{k}.$
For any $a\in \mathcal{F}_{3}\cap A_{m}$(since $w\in
\mathcal{F}_{3}$),$e_{1}w^{*}e_{1}(u_{1}^{k})^{*}e_{1}ae_{1}u_{1}^{k}e_{1}we_{1}\approx_{3\eta/k}e_{1}ae_{1}.$
 By (2),(3) above, it then follows from Lemma
2.9,  there is a monomorphism $\varphi: A_{m}\rightarrow
e_{1}Ae_{1}$ such that $\|\varphi(a)-e_{1}ae_{1}\|<\delta\|a\|$ for
all $a\in A_{m}.$

By applying Lemma 2.9 of \cite{H.Lin3}, we obtain unitaries
$x,x_{1},x_{2},\cdots,x_{k-1}\in U_{0}(e_{1}Ae_{1})$ such that
$\|x-e_{1}we_{1}\|<\delta,$ $\|x_{i}-e_{1}w_{i}e_{1}\|<\delta,$
$x=x_{1}x_{2}\cdots x_{k-1}$ and $x^{*}(u_{1}^{k})^{*}au_{1}^{k}x=a$
for all $a\in \varphi (A_{m}).$

 Let
 $Z=\sum_{i=1}^{k}e_{i}u_{1}^{k+1-i}x_{i}(u_{1}^{k-i})^*+(1-p)u_{1}.$
 Define $B=\varphi(A_{m}),$ then
 \begin{eqnarray*}\|Z-u_{1}\|&\leq &\max_{i}\{\|x_{i}-e_{1}\|\}\leq \max_{i}\{\|x_{i}-e_{1}w_{i}e_{1}\|+\|e_{1}w_{i}e_{1}-e_{1}\|\}\\ &<&\delta+\eta/k+\pi
 /(k-2),\end{eqnarray*} $
 (Z^{k})^{*}bZ^{k}=b$ for all $b\in B$ and
 $$(Z^{i})^{*}e_{1}Z^{i}\leq
 e_{i+1}, Z^{i}=u_{1}^{k}(x_{1}x_{2}\cdots x_{i})(u_{1}^{k-i})^{*},
 i=1,2,\cdots,k(e_{k+1}=e_{1}).$$

 Write $B=C_{1}\oplus C_{2}\oplus \cdots C_{N}$ and
 $\{c_{is}^{(j)}\}$ be the matrix units for
 $C_{j},j=1,2,\cdots,N,$ where $C_{j}=M_{R(j)}$ and put $q=1_{B}.$

 Define $D_{0}=B\oplus\oplus_{i=1}^{k-1}Z^{i*}BZ^{i},$
 and $D_{1}$ the $C$*-subalgebra generated by $B$ and $c_{ss}^{(j)}Z^{i}, s=1,2,\cdots,R(j), j=1,2,\cdots,N$
and $i=0,1,2,\cdots,k-1.$ Then $D_{1}\cong B\otimes M_{k}$ and
$D_{1}\supset D_{0}.$

 Define $q_{ss}^{(j)}=\sum_{i=0}^{k-1}Z^{i*}c_{ss}^{(j)}Z^{i},$ $q^{(j)}=\sum_{s=1}^{R(j)}q_{ss}^{(j)}$ and
 $Q=\sum_{j=1}^{N}q^{(j)}=1_{D_{1}}.$ Note that
 $Q=\sum_{i=0}^{k-1}(Z^{i})^{*}qZ^{i}.$ Note that
\begin{eqnarray*}q_{ss}^{(j)}Z&=&(\sum_{i=0}^{k-1}Z^{i*}c_{ss}^{(j)}Z^{i})Z=Z\sum_{i=0}^{k-1}(Z^{i+1})^{*}c_{ss}^{(j)}Z^{i+1}\\&=&
 Z(\sum_{i=1}^{k-1}Z^{i*}c_{ss}^{(j)}Z^{i}+c_{ss}^{(j)})=Zq_{ss}^{(j)}.\end{eqnarray*}

  It follows from Lemma 2.11 of \cite{H.Lin3} that
  $c_{11}^{(j)},c_{11}^{(j)}Z^{i}$ and
  $c_{11}^{(j)}Z^{k}c_{11}^{(j)}$ generate a $C$*-subalgebra which
  is isomorphic to $C(X_{j})\otimes M_{k}$ for some compact subset
  $X_{j}\subset S^{1}.$ Moreover, $q_{ss}^{(j)}Zq_{ss}^{(j)}$ is
  in the $C$*-subalgebra. Let $D$ be the $C$*-subalgebra generated
  by $D_{1}$ and $c_{11}^{(j)}Z^{k}c_{11}^{(j)}.$ Then $D\cong \oplus_{j=1}^{N}C(X_{j})\otimes B\otimes
  M_{k}.$ It  follows  that $q^{(j)}$ and $Q$ commutes
  with $Z$. Therefore $QZQ\in D.$ Thus,
\begin{eqnarray*}
  \|Qu-uQ\|&\leq & \|Qu-Qu_{1}\|+\|Qu_{1}-QZ\|+\|ZQ-u_{1}Q\|+\|u_{1}Q-uQ\|\\
  &<&4\eta/k+2\delta+2\pi/(k-2)<\varepsilon.
  \end{eqnarray*}

From $QZQ\in D,$ we also have $QuQ\in _{\varepsilon} D.$

For $b\in \mathcal{F}_{0},$ we compute that {\small
\begin{eqnarray*}&&(Z^{i})^{*}q(Z^{i})b=(Z^{i})^{*}qu_{1}^{k}(x_{1}x_{2}\cdots x_{i})(u_{1}^{k-i})^{*}b\\&\approx& _{k\delta+2\eta}(Z^{i})^{*}qu_{1}^{k}(w_{1}w_{2}\cdots
  w_{i})(u^{k-i})^{*}b\\
  &=& (Z^{i})^{*}qu_{1}^{k}(w_{1}w_{2}\cdots w_{i})
  (u^{k-i})^{*}bu^{k-i}(w_{1}w_{2}\cdots w_{i})^{*}(u_{1}^{k})^{*}[u^{k-i}(w_{1}w_{2}\cdots w_{i})^{*}
  (u_{1}^{k})^{*}]^{*}
  \end{eqnarray*}}

  Put $c_{i}=(u^{k-i})^{*}bu^{k-i},$ then $c_{i}\in
  \mathcal{F}_{1}.$ There is $a_{i}\in \mathcal{G}_{0}\subset
  A_{m}$ such that $\|c_{i}-a_{i}\|<\delta.$

Since $(w_{1}w_{2}\cdots w_{i})\mathcal{F}_{1}(w_{1}w_{2}\cdots
w_{i})^{*}\subset \mathcal{F}_{3},$ then {\small
\begin{eqnarray*}&
&(Z^{i})^{*}qu_{1}^{k}(w_{1}w_{2}\cdots
w_{i})(u^{k-i})^{*}bu^{k-i}(w_{1}w_{2}\cdots
w_{i})^{*}(u_{1}^{k})^{*}[u^{k-i}(w_{1}w_{2}\cdots
w_{i})^{*}(u_{1}^{k})^{*}]^{*}\\
&=&(Z^{i})^{*}qu_{1}^{k}(w_{1}w_{2}\cdots
w_{i})c_{i}(w_{1}w_{2}\cdots
w_{i})^{*}(u_{1}^{k})^{*}[u^{k-i}(w_{1}w_{2}\cdots
w_{i})^{*}(u_{1}^{k})^{*}]^{*}\\
&\approx&_{\delta} (Z^{i})^{*}qu_{1}^{k}(w_{1}w_{2}\cdots
w_{i})a_{i}(w_{1}w_{2}\cdots
w_{i})^{*}(u_{1}^{k})^{*}[u^{k-i}(w_{1}w_{2}\cdots
w_{i})^{*}(u_{1}^{k})^{*}]^{*}\\
&\approx&_{\delta} (Z^{i})^{*}e_{1}u_{1}^{k}(w_{1}w_{2}\cdots
w_{i})a_{i}(w_{1}w_{2}\cdots
w_{i})^{*}(u_{1}^{k})^{*}[u^{k-i}(w_{1}w_{2}\cdots
w_{i})^{*}(u_{1}^{k})^{*}]^{*}\\
&\approx&_{\eta/k} (Z^{i})^{*}u_{1}^{k}(w_{1}w_{2}\cdots
w_{i})a_{i}(w_{1}w_{2}\cdots
w_{i})^{*}(u_{1}^{k})^{*}e_{1}[u^{k-i}(w_{1}w_{2}\cdots
w_{i})^{*}(u_{1}^{k})^{*}]^{*}\\
&\approx&_{\delta} (Z^{i})^{*}u_{1}^{k}(w_{1}w_{2}\cdots
w_{i})c_{i}(w_{1}w_{2}\cdots
w_{i})^{*}(u_{1}^{k})^{*}q[u^{k-i}(w_{1}w_{2}\cdots
w_{i})^{*}(u_{1}^{k})^{*}]^{*}\\
&\approx&_{\delta} (Z^{i})^{*}u_{1}^{k}(w_{1}w_{2}\cdots
w_{i})(u^{k-i})^{*}bu^{k-i}(w_{1}w_{2}\cdots
w_{i})^{*}(u_{1}^{k})^{*}q[u^{k-i}(w_{1}w_{2}\cdots
w_{i})^{*}(u_{1}^{k})^{*}]^{*}\\
&\approx&_{k\delta+2\eta} b(Z^{i})^{*}qZ^{i}.\end{eqnarray*}}

 Hence
 $$\|(Z^{i})^{*}qZ^{i}b-b(Z^{i})^{*}qZ^{i}\|<2(k\delta+2\eta+\delta+\delta)+\eta/k<\varepsilon/k,k=0,1,\cdots,k-1.$$
 Therefore, for $b\in \mathcal{F}_{0},$ $\|Qb-bQ\|<k\cdot(\varepsilon/k)=\varepsilon.$

 It follows that $\|Qa-aQ\|<\varepsilon$ for all $a\in \mathcal{F}.$

 For any $b\in \mathcal{F}_{0},$ a same estimation above shows
 that

 $$\|qZ^{i}b(Z^{i})^{*}q-qu_{1}^{k}(w_{1}w_{2}\cdots w_{i})(u^{k-i})^{*}bu^{k-i}(w_{1}w_{2}\cdots
 w_{i})^{*}(u_{1}^{k})^{*}q\|<2k\delta+4\eta$$

 However,
 $qu_{1}^{k}(w_{1}w_{2}\cdots w_{i})(u^{k-i})^{*}bu^{k-i}(w_{1}w_{2}\cdots w_{i})^{*}(u_{1}^{k})^{*}q\in_{\delta+2\delta+4\eta/k}B.$

 It follows that, for $b\in \mathcal{F}_{0},$

 $(Z^{i})^{*}qZ^{i}b(Z^{i})^{*}qZ^{i}\in _{\varepsilon/k}(Z^{i})^{*}BZ^{i}, i=0,1,\cdots,k-1.$

 we obtain that $QbQ\in_{\varepsilon} D_{1}\subset D.$

 we obtain that $QaQ\in_{\varepsilon}D$ for all $a\in \mathcal{F}.$

Because $[1-\sum_{i=1}^{k}e_{i}]=[1-p]\leq_{\alpha} [r_{0}]$ in $A$,
there exist the mutually orthogonal projections $p_{i},$ the
mutually orthogonal positive elements $a_{i}$ and $s_{i}\in
\mathbb{Z}$ for $i=1,2,\cdots,n$ such that
$p=\sum_{i=1}^{n}p_{i},\{a_{i}\}_{i=1}^{n}$ belong to the hereditary
$C$*-subalgebra generated by $r_0$,  and
$[\alpha^{s_{i}}(p_{i})]\leq [a_{i}],$ $i=1\cdots n.$
 Because
$[\alpha^{s_{i}}(p_{i})]=[u^{s_{i}}p_{i}(u^{s_{i}})^{*}]=[p_{i}]$ in
$A\rtimes_{\alpha}\mathbb{Z}$. we obtain that
$[1-\sum_{i=1}^{k}e_{i}]\leq [r_{0}]$ in
$A\rtimes_{\alpha}\mathbb{Z}$

 We can compute that

 $[1-Q]\leq [1-\sum_{i=1}^{k}e_{i}]\leq [r_{0}]\leq[b_{0}].$

 So TR$(A\rtimes_{\alpha}\mathbb{Z})\leq 1.$

The order on projection over $A\rtimes_{\alpha}\mathbb{Z}$ is
determined by traces by Theorem 3.7.2 of \cite{H.Lin2}.

By applying Theorem 3.2, we have
RR$(A\rtimes_{\alpha}\mathbb{Z})=0$.
 By Lemma 3.2 of \cite{H.Lin3}, we conclude that TR$(A\rtimes_{\alpha}\mathbb{Z})=0.$
\end{proof}

\begin{Cor} Let $A$ be a unital AF-algebra, Suppose that
$\alpha\in$Aut$(A)$ has the tracial cyclic Rokhlin property and $A$
is $\alpha$-simple. Suppose also that there is an integer $J\geq1$
such that $\alpha^{J}_{*0}=$id$_{K_{0}(A)}$. Then the restriction
map is a bijection from the tracial states of
$A\rtimes_{\alpha}\mathbb{Z}$ to the $\alpha$-invariant tracial
states of $A$.
\end{Cor}

\begin{proof} Since $A$ has real rank zero and
$A\rtimes_{\alpha}\mathbb{Z}$ also has real rank zero by Theorem
3.3, this follows from Proposition 2.2 of \cite{A.Kishimoto3}.

\end{proof}

\begin{Example}Let $A=A_{0}\oplus A_{0},$ where $A_{0}$ is an
 infinite dimensional unital simple AF-algebra. Let $\beta\in
 $Aut$(A_{0})$ be an approximately inner automorphism of $A_{0}$ and have the traical cyclic Rohklin property in
 \cite{H.Lin3}.  Define $\alpha\in $ Aut$(A)$ by
 $\alpha(a,b)=(\beta(b),\beta(a)),$ then TR$(A\rtimes_{\alpha}\mathbb{Z})=0$.\end{Example}

Obviously, $A$ is $\alpha$-simple. Because $\beta$ is an
approximately inner automorphism of $A_{0}$, therefore
$\beta_{*0}=$id$_{K_{0}(A_{0})},$ then we have
$\alpha^{2}_{*0}=$id$_{K_{0}(A)}.$

Because $\beta$ is an approximately inner automorphism of $A_{0}$
and has the traical cyclic Rohklin property in
 \cite{H.Lin3}, furthermore by applying Lemma 2.8 of \cite{H.Lin3}, it is easy to verify that  $\alpha$ has the traical cyclic Rohklin
 property in this paper.

So $(A,\alpha)$
 satisfies the conditions of Theorem 3.3, then we have
 TR$(A\rtimes_{\alpha}\mathbb{Z})=0$.

\

\

{\scshape Department of Mathematics, East China Normal University,
Shanghai 200241, P.R.CHINA.}

{\it E-mail address}: huajiajie2006@hotmail.com
\end{document}